\documentclass[10pt]{article} \usepackage{amsfonts} \usepackage{amssymb}\usepackage{amsmath}\usepackage{amsfonts}

\makeatletter
\renewcommand{\@makefnmark}{}
\makeatother
\begin{document}
\baselineskip=10pt
\pagestyle{plain}
{\Large

\medskip
\medskip
\medskip

\footnote{

\hspace{2mm}Keywords: Dirac operator, degenerate boundary conditions, spectrum}

\centerline {\bf On the spectrum of  $2\times 2$ Dirac operator}
\centerline {\bf   with degenerate boundary conditions}
\medskip
\medskip
\medskip
\medskip
\medskip

\centerline { Alexander Makin}
\medskip
\medskip
\medskip
{\normalsize
\centerline {Peoples Friendship University of Russia}
 \centerline {	117198, Miklukho-Maklaya str. 6, Moscow, Russia}

\medskip
\medskip
\medskip

\medskip
\medskip
\begin{quote}{\normalsize
We study the spectral problem for the Dirac operator with degenerate boundary conditions and a complex-valued summable potential. Sufficient conditions are found under which the spectrum of the problem under consideration coincides with the spectrum of the corresponding unperturbed operator.}
   \end{quote}

\centerline {\bf 1. Introduction}

\medskip
\medskip
\medskip
\medskip

In the  present paper, we study the Dirac system
$$
  B\mathbf{y}'+ V\mathbf{y}=\lambda \mathbf{y},\eqno(1)
 $$
 where
 $$
  B=\begin{pmatrix}
 0&1\\
 -1&0
 \end{pmatrix},\quad  V=\begin{pmatrix}
 p(x)&q(x)\\
 q(x)&-p(x)
 \end{pmatrix},
$$
the functions $p, q\in L_1(0,\pi)$,   with two-point boundary conditions

$$
\begin{array}{c}
U_1(\mathbf{y})= a_{11}y_1(0)+a_{12}y_2(0)+a_{13}y_1(\pi)+ a_{14}y_2(\pi)=0,\\ U_2(\mathbf{y})= a_{21}y_1(0)+a_{22}y_2(0)+a_{23}y_1(\pi)+ a_{24}y_2(\pi)=0,
\end{array}\eqno(2)
$$
where
the coefficients $a_{jk}$ are arbitrary complex numbers,
and rows of  the matrix
$$
A=\begin{pmatrix}
 a_{11}&a_{12}&a_{13} &a_{14}\\
 a_{21}&a_{22}&a_{23} &a_{24}
\end{pmatrix}
$$
are linearly independent.

We consider the operator $\mathbb{L}\mathbf{y}=B\mathbf{y}'+V\mathbf{y}$ as a linear operator on the space $\mathbb{H}=L_2(0,\pi)\oplus L_2(0,\pi)$,
with the domain $D(\mathbb{L})=\{\mathbf{y}\in W_1^1[0,\pi]:\, \mathbb{L}\mathbf{y}\in \mathbb{H}$, $U_j(\mathbf{y})=0$ $(j=1,2)\}$.

Denote  by $J_{jk}$ the determinant composed of the $j$th and $k$th columns  of the matrix $A$.
Denote $J_0=J_{12}+J_{34}$, $J_1=J_{14}-J_{23}$, $J_2=J_{13}+J_{24}$.

Boundary conditions (2) are called degenerate if

$$
J_1=J_2=0; \quad J_0=0,\quad J_1+iJ_2\ne0, \quad J_1-iJ_2=0;\quad J_0=0, \quad J_1+iJ_2=0, \quad J_1-iJ_2\ne0,
$$
otherwise they are  nondegenerate.

There is an enormous literature related to the spectral theory for Dirac operators with nondegenerate boundary conditions. The case of degenerate conditions has been  investigated much less although in the last decade interest in the study of these spectral problems has increased significantly [see 1-9 and the references therein].
The main goal of present paper is to establish conditions on the potential $V$ under which the spectrum of the problem under consideration with degenerate boundary conditions coincides with the spectrum of the corresponding unperturbed operator (3)
$$
B\mathbf{y}'=\lambda\mathbf{y},\quad U(\mathbf{y})=0.\eqno(3)
$$

\medskip
\medskip
\medskip
\centerline {\bf 2. Main results}
\medskip
\medskip
\medskip

{\bf Theorem 1.} {\it Suppose the following conditions are valid
$$
J_{14}=J_{23}=J_{13}+J_{24}=0,\eqno(4)
$$

$$
p(\pi-x)=-p(x), \quad q(\pi-x)=q(x), \eqno(5)
$$
where $0\le x\le\pi$. Then the spectrum of problem (1) (2) coincides with the spectrum of the corresponding unperturbed operator (3).}

Proof.
First of all, we rewrite system (1) in scalar form
$$
\left\{
\begin{array}{rcl}
y_2'+p(x)y_1+q(x)y_2=\lambda y_1\\
-y_1'+q(x)y_1-p(x)y_2=\lambda y_2.\\
\end{array}
\right.\eqno(6)
$$
Denote by
 \[
E(x,\lambda)=\begin{pmatrix}
e_{11}(x,\lambda)&e_{12}(x,\lambda)\\
e_{21}(x,\lambda)&e_{22}(x,\lambda)
\end{pmatrix}\eqno(7)
\]
the matrix of the fundamental solution system  to  system
(1)
with boundary condition
$
E(\frac{\pi}{2},\lambda)=I
$, where $I$ is the unit matrix.
It is well known that
$$
e_{11}(x,\lambda)e_{22}(x,\lambda)-e_{12}(x,\lambda)e_{21}(x,\lambda)=1\eqno(8)
$$
for any $x,\lambda$.

Substituting the first column of matrix (7) in system (6), we obtain
$$
\left\{
\begin{array}{rcl}
e_{21}'(x,\lambda)+p(x)e_{11}(x,\lambda)+q(x)e_{21}(x,\lambda)=\lambda e_{11}(x,\lambda)\\
-e_{11}'(x,\lambda)+q(x)e_{11}(x,\lambda)-p(x)e_{21}(x,\lambda)=\lambda e_{21}(x,\lambda).\\
\end{array}
\right.\eqno(9)
$$
Replacing $x=\pi-t$ in relations (9), we find
$$
\left\{
\begin{array}{rcl}
-e_{21}'(\pi-t,\lambda)+p(\pi-t)e_{11}(\pi-t,\lambda)+q(\pi-t)e_{21}(\pi-t,\lambda)=\lambda e_{11}(\pi-t,\lambda)\\
e_{11}'(\pi-t,\lambda)+q(\pi-t)e_{11}(\pi-t,\lambda)-p(\pi-t)e_{21}(\pi-t,\lambda)=\lambda e_{21}(\pi-t,\lambda).\\
\end{array}
\right.\eqno(10)
$$
It follows from (5) and (10) that
$$
\left\{
\begin{array}{rcl}
-e_{21}'(\pi-t,\lambda)-p(t)e_{11}(\pi-t,\lambda)+q(t)e_{21}(\pi-t,\lambda)=\lambda e_{11}(\pi-t,\lambda)\\
e_{11}'(\pi-t,\lambda)+q(t)e_{11}(\pi-t,\lambda)+p(t)e_{21}(\pi-t,\lambda)=\lambda e_{21}(\pi-t,\lambda).\\
\end{array}
\right.\eqno(11)
$$
Denote $z_2(t,\lambda)=e_{11}(\pi-t,\lambda)$, $z_1(t,\lambda)=e_{21}(\pi-t,\lambda)$. It follows from (11) that

$$
\left\{
\begin{array}{rcl}
z_2'+p(t)z_1+q(t)z_2=\lambda z_1\\
-z_1'+q(t)z_1-p(t)z_2=\lambda z_2.\\
\end{array}
\right.\eqno(12)
$$
Obviously, $z_2(\frac{\pi}{2},\lambda)=1$, $z_1(\frac{\pi}{2},\lambda)=0$,  therefore, by virtue of the uniqueness of the solution to  the Cauchy problems for systems (9) and (12) $z_2(t,\lambda)=e_{22}(t,\lambda)$, $z_1(t,\lambda)=e_{21}(t,\lambda)$, hence, we obtain
$$
e_{22}(t,\lambda)=e_{11}(\pi-t,\lambda), \quad e_{21}(t,\lambda)=e_{21}(\pi-t,\lambda).
$$
The last relations imply
$$
e_{21}(\pi)=e_{12}(0),\quad e_{12}(\pi)=e_{21}(0),\quad e_{11}(\pi)=e_{22}(0),\quad e_{22}(\pi)=e_{11}(0).\eqno(13)
$$

The eigenvalues of problem (1), (2) are the roots of the characteristic equation
$$
\Delta(\lambda)=0,
$$
where
$$
\Delta(\lambda)=
\left|\begin{array}{cccc}
U_1(E^{[1]}(\cdot,\lambda))&U_1(E^{[2]}(\cdot,\lambda))\\
U_2(E^{[1]}(\cdot,\lambda))&U_2(E^{[2]}(\cdot,\lambda))\\
\end{array}
\right|,
$$
$E^{[k]}(x,\lambda)$ is the $k$th column of matrix (7).
Simple computations together with relations (4), (8), and (13) show that
$$
\Delta(\lambda)=
\left|\begin{array}{cccc}
a_{11}e_{11}(0)+a_{12}e_{21}(0)+a_{13}e_{11}(\pi)+ a_{14}e_{21}(\pi)&a_{11}e_{12}(0)+a_{12}e_{22}(0)+a_{13}e_{12}(\pi)+ a_{14}e_{22}(\pi)\\
a_{21}e_{11}(0)+a_{22}e_{21}(0)+a_{23}e_{11}(\pi)+ a_{24}e_{21}(\pi)&a_{21}e_{12}(0)+a_{22}e_{22}(0)+a_{23}e_{12}(\pi)+ a_{24}e_{22}(\pi)\\
\end{array}
\right|
$$
$$
\begin{array}{c}
=[a_{11}e_{11}(0)+a_{12}e_{21}(0)+a_{13}e_{11}(\pi)+ a_{14}e_{21}(\pi)][a_{21}e_{12}(0)+a_{22}e_{22}(0)+a_{23}e_{12}(\pi)+ a_{24}e_{22}(\pi)]\\-
[a_{21}e_{11}(0)+a_{22}e_{21}(0)+a_{23}e_{11}(\pi)+ a_{24}e_{21}(\pi)][a_{11}e_{12}(0)+a_{12}e_{22}(0)+a_{13}e_{12}(\pi)+ a_{14}e_{22}(\pi)]\\\\

=[a_{11}e_{11}(0)+a_{12}e_{21}(0)+a_{13}e_{22}(0)+ a_{14}e_{12}(0)][a_{21}e_{12}(0)+a_{22}e_{22}(0)+a_{23}e_{21}(0)+ a_{24}e_{11}(0)]\\-
[a_{21}e_{11}(0)+a_{22}e_{21}(0)+a_{23}e_{22}(0)+ a_{24}e_{12}(0)][a_{11}e_{12}(0)+a_{12}e_{22}(0)+a_{13}e_{21}(0)+ a_{14}e_{11}(0)]\\\\=

[a_{11}e_{11}(0)a_{21}e_{12}(0)+a_{11}e_{11}(0)a_{22}e_{22}(0)+a_{11}e_{11}(0)a_{23}e_{21}(0)+a_{11}e_{11}(0) a_{24}e_{11}(0)\\+
a_{12}e_{21}(0)a_{21}e_{12}(0)+a_{12}e_{21}(0)a_{22}e_{22}(0)+a_{12}e_{21}(0)a_{23}e_{21}(0)+a_{12}e_{21}(0)a_{24}e_{11}(0)\\+
a_{13}e_{22}(0)a_{21}e_{12}(0)+a_{13}e_{22}(0)a_{22}e_{22}(0)+a_{13}e_{22}(0)a_{23}e_{21}(0)+a_{13}e_{22}(0) a_{24}e_{11}(0)\\+
a_{14}e_{12}(0)a_{21}e_{12}(0)+a_{14}e_{12}(0)a_{22}e_{22}(0)+a_{14}e_{12}(0)a_{23}e_{21}(0)+a_{14}e_{12}(0)a_{24}e_{11}(0)]\\-
[a_{21}e_{11}(0)a_{11}e_{12}(0)+a_{21}e_{11}(0)a_{12}e_{22}(0)+a_{21}e_{11}(0)a_{13}e_{21}(0)+a_{21}e_{11}(0)a_{14}e_{11}(0)\\+
a_{22}e_{21}(0)a_{11}e_{12}(0)+a_{22}e_{21}(0)a_{12}e_{22}(0)+a_{22}e_{21}(0)a_{13}e_{21}(0)+a_{22}e_{21}(0)a_{14}e_{11}(0)\\+
a_{23}e_{22}(0)a_{11}e_{12}(0)+a_{23}e_{22}(0)a_{12}e_{22}(0)+a_{23}e_{22}(0)a_{13}e_{21}(0)+a_{23}e_{22}(0)a_{14}e_{11}(0)\\+
a_{24}e_{12}(0)a_{11}e_{12}(0)+a_{24}e_{12}(0)a_{12}e_{22}(0)+a_{24}e_{12}(0)a_{13}e_{21}(0)+a_{24}e_{12}(0)a_{14}e_{11}(0)]\\\\=

e_{11}(0)e_{22}(0)(a_{11}a_{22}+a_{13} a_{24}-a_{21}a_{12}-a_{23}a_{14})+
e^2_{11}(0)(a_{11}a_{24}-a_{21}a_{14})+e^2_{22}(0)(a_{13}a_{22}-a_{23}a_{12})\\+
e_{12}(0)e_{21}(0)(a_{12}a_{21}+a_{14}a_{23}-a_{24}a_{13}-a_{22}a_{11})+e_{11}(0)e_{12}(0)(a_{11}a_{21}+a_{14}a_{24}-a_{21}a_{11}-a_{24}a_{14})\\+
e_{11}(0)e_{21}(0)(a_{11}a_{23}+a_{12}a_{24}-a_{21}a_{13}-a_{22}a_{14})\\+
e^2_{21}(0)(a_{12}a_{23}-a_{22}a_{13})+e^2_{12}(0)(a_{14}a_{21}-a_{24}a_{11})\\+
e_{21}(0)e_{22}(0)(a_{12}a_{22}+a_{13}a_{23}-a_{22}a_{12}-a_{23}a_{13})\\+
e_{22}(0)e_{12}(0)(a_{13}a_{21}+a_{14}a_{22}-a_{23}a_{11}-a_{24}a_{12})\\\\=

e_{11}(0)e_{22}(0)(J_{12}+J_{34})+e^2_{11}(0)J_{14}+e^2_{22}(0)J_{32}\\+
e_{12}(0)e_{21}(0)(J_{21}+J_{43})+e_{11}(0)e_{21}(0)(J_{13}+J_{24})\\+
e^2_{21}(0)J_{23}+e^2_{12}(0)J_{41}+e_{22}(0)e_{12}(0)(J_{31}+J_{42})\\\\=

[e_{11}(0)e_{22}(0)-e_{12}(0)e_{21}(0)](J_{12}+J_{34})+e^2_{11}(0)J_{14}+e^2_{22}(0)J_{32}+e^2_{21}(0)J_{23}+e^2_{12}(0)J_{41}\\+
[e_{11}(0)e_{21}(0)-e_{22}(0)e_{12}(0)](J_{13}+J_{24)})\\\\=

J_{12}+J_{34}+[e^2_{11}(0)-e^2_{12}(0)]J_{14}+[e^2_{21}(0)-e^2_{22}(0)]J_{23}+
[e_{11}(0)e_{21}(0)-e_{22}(0)e_{12}(0)](J_{13}+J_{24})=J_{12}+J_{34}.

\end{array}
$$
This completes the proof.
\medskip
\medskip

{\bf Remark 1.} Condition (4) holds if, for example,

$$
A=\begin{pmatrix}
 1&0&0 &b\\
 0&1&b &0
\end{pmatrix}
 \quad\mbox{or}\quad
 A=\begin{pmatrix}
 0&b&1 &0\\
 b&0&0 &1
\end{pmatrix}.\eqno(14)
$$
In both cases $J_{12}+J_{34}=1-b^2$, hence, if $b^2\ne1$ the spectrum is empty, and if $b^2=1$ the spectrum fills all complex plane.
Notice, that if $b=0$ conditions (14) are the Cauchy boundary conditions. The  Cauchy problem has no spectrum for any potential $V$.

{\bf Remark 2.} Let us consider system (1) with boundary  conditions

$$
\tilde A=\begin{pmatrix}
 1&0&0 &1\\
 0&1&1 &0
\end{pmatrix}.\eqno(15)
$$
Denote by
$$
\tilde E(x,\lambda)=\begin{pmatrix}
c_1(x,\lambda)&-s_2(x,\lambda)\\
s_1(x,\lambda)&c_2(x,\lambda)
\end{pmatrix}
$$
the matrix of the fundamental solution system  to  system (1)
with boundary condition
$
\tilde E(0,\lambda)=I
$, where
$I$ is the unit matrix. It is easy to see that $J_{14}=J_{23}=J_{12}+J_{34}=0, J_{13}=1, J_{24}=-1$.
Trivial computation shows that  the characteristic determinant of problem (1), (15) can be reduced to the form
$$
\tilde\Delta(\lambda)=s_1(\pi,\lambda)-s_2(\pi,\lambda).
$$
Let condition (5) hold. From Theorem 1, we get that any complex $\lambda$ is an eigenvalue of problem (1), (15), hence, $\tilde\Delta(\lambda)\equiv0$, therefore, $s_1(\pi,\lambda)\equiv s_2(\pi,\lambda)$.

\medskip
\medskip
\medskip

email: alexmakin@yandex.ru

\end{document}